\documentclass[11pt]{article}
\pdfoutput=1 

\usepackage[utf8]{inputenc}
\RequirePackage[displaymath, mathlines, pagewise]{lineno}

\usepackage{algorithm, algorithmicx, algpseudocode}
\usepackage{listings}%
\usepackage{alltt}
\usepackage{caption}
\usepackage[margin=1in]{geometry}
\usepackage{amsfonts}
\usepackage{amsmath}
\usepackage{color}
\usepackage{cite}

\RequirePackage{graphicx}

\DeclareMathOperator{\sgn}{sgn}
\DeclareMathOperator*{\argmax}{arg\,max}
\DeclareMathOperator*{\argmin}{arg\,min}
\newcommand{\change}[1]{{\color{black} #1}}

\begin{document}

\title{A Continuous Time Formulation of Stochastic Dual Control To Avoid the Curse of Dimensionality}

\author{Martin P\'{e}ron$^{1,2}$*, Christopher M. Baker$^{1,2,3,4}$, Barry D. Hughes$^4$ and Iadine Chad\`{e}s$^2$\\[\baselineskip]
\it $^1$School of Mathematical Sciences, Queensland University of Technology,\\
\it Brisbane, Queensland 4000, Australia\\
\it $^2$ Land and Water, CSIRO, Ecosciences Precinct,\\
\it  Dutton Park, Queensland 4102, Australia\\
\it $^3$ School of Biological Sciences, University of Queensland,\\
\it  Saint Lucia, Queensland 4072, Australia\\
\it $^4$ School of Mathematics and Statistics, University of Melbourne,\\
 Parkville, Victoria 3010, Australia\\[\baselineskip]
  *email {\tt martin.peron@laposte.net}
}

\maketitle

\begin{abstract}
\noindent Dual control denotes a class of control problems where the parameters governing the system are imperfectly known. The challenge is to find the optimal balance between probing, i.e. exciting the system to understand it more, and caution, i.e. selecting conservative controls based on current knowledge to achieve the control objective. Dynamic programming techniques can achieve this optimal trade-off. However, while dynamic programming performs well with discrete state and time, it is not well-suited to problems with continuous time-frames or continuous or unbounded state spaces. Another limitation is that multidimensional states often cause the dynamic programming approaches to be intractable. In this paper, we investigate whether continuous-time optimal control tools could help circumvent these caveats whilst still achieving the probing--caution balance. We introduce a stylized problem where the state is governed by one of two differential equations. It is initially unknown which differential equation governs the system, so we must simultaneously determine the `true' differential equation and control the system to the desired state. We show how this problem can be transformed to apply optimal control tools, and compare the performance of this approach to a dynamic programming approach. Our results suggest that the optimal control algorithm rivals dynamic programming on small problems, achieving the right balance between aggressive and smoothly varying controls. In contrast to dynamic programming, the optimal control approach remains tractable when several states are to be controlled simultaneously.
\end{abstract}

\begin{center}
stochastic differential equations, dynamic programming, adaptive control
\end{center}

\thispagestyle{empty}

\newpage

\section{Introduction}
Many real-world control problems involve uncertainty, caused by a too complex or poorly understood system or inaccurate measurements. So, decision makers are often uncertain about how their actions affect the system \cite{astrom_adaptive_2008}. This uncertainty must be accounted for to make optimal decisions. In the control literature, where the aim is typically to control a system towards a desired state, this problem is called dual control \cite{astrom_adaptive_2008, bertsekas_dynamic_1995, feldbaum_optimal_1965}, while in environmental sciences it is called adaptive management \cite{walters_ecological_1978,chades_optimization_2017}. The problem arises in a broad range of fields, including industry \cite{astrom_theory_1983}, conservation \cite{chades_momdps_2012,runge_active_2013} and natural resource management \cite{johnson_conditions_2002,frederick_choosing_1995}. In both the optimal control and environmental science communities, the problem is modeled by uncertain parameters augmenting the system state. Learning these parameters is usually not the primary control objective, rather, the objective usually only depends on the state and control. Thus, the optimal solution trades off decisions to better understand the system and decisions to guide the system towards a better state. Hence, \textit{informative} controls (i.e. improving knowledge) should only be chosen over more \textit{rewarding} controls if the long-term benefits of learning outweigh the potential short-term loss of performance \cite{walters_ecological_1978}. In dual control terminology, this is  `finding the optimal balance between \textit{probing} and \textit{caution}' \cite{bertsekas_dynamic_1995}. 

Achieving the balance between probing and caution is not an easy task. Researchers have mostly focused on solving discrete-time and discrete-state problems \cite{astrom_adaptive_2008,bertsekas_dynamic_1995, chades_optimization_2017}, as they are deemed easier to solve. There exist continuous-time exceptions but with no attempt at finding the optimal trade-off \cite{naik_robust_1992}: the control follows the certainty equivalence principle, i.e. the control is chosen as if the current estimate of the uncertain parameters were true \cite{bertsekas_dynamic_1995}. In contrast, it is in theory possible to find the optimal discrete-time control by modeling the problem as a Markov decision problem (or a variant) and solving it by using stochastic dynamic programming \cite{astrom_adaptive_2008,chades_momdps_2012,puterman_markov_1994}.

Unfortunately, Markovian frameworks usually require specifying a finite set of states, which makes it challenging to implement such frameworks when the state is unbounded or multidimensional (curse of dimensionality) \cite{bellman_dynamic_1957}. This is further reinforced by the PSPACE-complete complexity of such problems \cite{chades_momdps_2012}. For these reasons, the optimal dual controller has even been said to be `impossible' to calculate for real-world processes \cite{astrom_adaptive_2008}.  

Many approaches have been explored to circumvent the curse of dimensionality. Methods drawing from approximate dynamic programming \cite{lee_approximate_2009}, or using hand-made control criteria in an attempt to balance probing and caution \cite{dumont_wood_1988}, are promising but they are approximate approaches with no performance guarantees. Another body of work revolves around model predictive control (MPC). MPC uses previous system states, inputs and outputs to predict good control actions over a finite receding horizon. Variants of MPC combine it with ideas from optimal experimental design \cite{bavdekar_stochastic_2016} or (in self-reflective MPC) incorporate system noise and measurement errors \cite{feng_real-time_2018,houska_self-reflective_2017}.
Although MPC is an attractive heuristic approach to dual control problems \cite{camacho_model_2013,thangavel_towards_2015,heirung_mpc_2013} and does entail some exploration of the state space, it does not achieve true optimization of the trajectory to a final state that lies beyond the finite horizon when the MPC process is begun.

Tools from continuous--time optimal control can help circumvent these caveats. Continuous--time optimal control problems can be solved by different methods, one of which is the Pontryagin minimum principle \cite{bertsekas_dynamic_1995}. This approach leads to differential equations that can be solved numerically to find the optimal control \cite{baker_placing_2016,hackbusch_numerical_1978,lenhart_optimal_2007}. Here we explore how continuous-time control tools could be applied to actively learning policies in a continuous-time, unbounded- and continuous-state setting. 

We introduce a simple continuous-time, unbounded- and continuous-state problem, with only two possible values for the uncertain parameter. Firstly, we show how to augment the \textit{physical state} with an unknown \textit{information state}, which represents the uncertain system dynamics. Secondly, we identify the stochastic differential equation that the information state satisfies. Thirdly, we show how to circumvent the stochasticity of our problem by replacing both the information and physical states by their expected values. Fourthly, we solve the resulting deterministic problem by an existing optimal control algorithm. Finally, we evaluate this approach through simulations and compare its performance to dynamic programming.

\section{Our simple dual-control problem vs a standard optimal control problem} 

In this section, we introduce the decision problem and the basic elements of solution methods. We then introduce our approach based on optimal control before comparing it to a dynamic programming approach.

In our model, the dynamics are governed by a diffusion process (in a stochastic differential equation framework). Our objective is to keep the state at zero. The challenge is that we are initially unsure whether the control affects the system state positively or negatively (parametric uncertainty). While this is a simple model, it is related to other control problems where parameters can change in sign, or the model itself is initially unknown \cite{dumont_wood_1988,bond_using_2009}. Also, this model is potentially useful for proof-of-concept testing for the applicability of continuous optimal control methods in this class of stochastic dual control problems.

\subsection{The problem: continuous-time dual control}

We aim to control a system with state $x(t) \in \mathbb{R}$ by choosing a certain control $u(t) \in [-U, U]$ for time $t \in [0,T]$. We use the concise notation $u(t)$ but the control may depend implicitly on the history of states and controls up to time $t$. The true state of the system, $x$, is governed by one of the two following stochastic differential equations: 
\begin{align} 
dx(t) &= u(t)dt + dB_t, \label{equation1} \\
dx(t) &= -u(t)dt + dB_t, \label{equation2} 
\end{align}
where $B_t$ is a Wiener Process, which satisfies for all $t,t' \in [0,T]$ with $t \leq t'$:
\begin{align}  
B_{t'}-B_t \sim \mathcal{N}(0, t'-t),
\end{align} 
independently of past values $B_s, s < t$. As usual, $\mathcal{N}(0,t)$ denotes the normal law with mean $0$ and variance $t$. In this problem, the state, $x(t)$, is perfectly observable. We assume that Equations \ref{equation1} and \ref{equation2} are equally likely to govern the dynamics of the state $x(t)$.

We set as our objective that the state $x(t)$ and the control $u(t)$ are both kept small in a suitable mean-square sense over the time interval $[0,T]$, conditioned on whether Equation \ref{equation1} or \ref{equation2} is true:
\begin{align}
\min_{u} \Bigl\{\dfrac{1}{2}\ \mathbf{E} \Big[
\int_0^T [x^2(t) + u^2(t)] dt \mid \mbox{Equation\ \ref{equation1} true} \Big] + \dfrac{1}{2}\ \mathbf{E} \Big[
\int_0^T [x^2(t) + u^2(t)] dt \mid \mbox{Equation\ \ref{equation2} true} \Big] \Bigr\}.
\label{equation3}
\end{align}
The effect of the control $u(t)$ on the state $x(t)$ is initially unknown and depends on which of Equation \ref{equation1} or \ref{equation2} is true. However, whether Equation \ref{equation1} or \ref{equation2} is true should become clearer over time, based on past observations of the variation of $x(t)$.

\subsection{Information state}

A common approach in both dual control (state augmentation \cite{bertsekas_dynamic_1995}) and adaptive management \cite{chades_optimization_2017}, is to create an \textit{information state}, $y_{true}$, representing the true underlying system dynamics. In our case, the state $y_{true}$ equals 1 if Equation \ref{equation1} is true and 0 if Equation \ref{equation2} is true. Since our knowledge is imperfect and may vary over time, we use the notation $y(t)$ to describe our belief that Equation \ref{equation1} is true at any time $t$. Note that $y_{true}$ is binary and should not be confused with $y(t)$, which is continuous within $[0,1]$. Our belief that Equation \ref{equation2} is true is $1-y(t)$. We will set the initial belief to $y({0})=1/2$ in the experiments to model the lack of prior information --- that is, assume the two state equations to be equally likely at $t=0$. 

In the next section we show how we can use tools from optimal control theory to solve this problem.

\subsection{Continuous-time optimal control}

A continuous-time optimal control problem \cite{bertsekas_dynamic_1995} aims to control, for time $t \in [0,T]$, a system in state $x(t) \in X$ by choosing a certain control $u(t)$. The objective is to minimize a certain cost $J$ defined as
\begin{align}  
J(u)=\int_0^T g(x(t), u(t))dt + f(x(T)).  \label{eq:OCobj1}
\end{align} 
The state equation is a differential equation that describes the evolution of the state $x(t)$ from the initial time $t=0$ to the final time $t=T$, under the action of the control $u(t)$: 
\begin{align}  
\frac{dx(t) }{dt} = h(x(t), u(t)).  \label{eq:state}
\end{align} 
The control may be bounded by a positive constant $U$:
\begin{align}  
|u(t)| \leq U.
\end{align} 
Formally, the objective is: 
\begin{align}  
\min_{u} J(u).  \label{eq:OCobj2}
\end{align}

The Pontryagin minimum principle is designed to solve the type of problem defined by Equations \ref{eq:OCobj1}--\ref{eq:OCobj2}. Details for this approach can be found in Appendix \ref{appendix:pontryagin}. Although the Pontryagin minimum principle is an elegant and efficient way to solve deterministic continuous-time control problems, it is not directly applicable to stochastic differential equations in general, including the problem we introduced. In the next section we show how to circumvent the stochasticity of this problem to enable the application of the Pontryagin minimum principle.

\section{Using tools from continuous-time optimal control}
 
Although the Pontryagin minimum principle cannot be applied directly, we shall show that a control strategy can still be based on it.

\subsection{Approach overview}\label{secApproachOverview}

First, we design an (approximate) deterministic optimal control problem which aims at capturing the uncertainty on the states. Although the resulting deterministic control profile can provide useful information, it is not sufficient in itself. The reason for this is the stochasticity of our original problem, which means an adaptive response is required. \change{To this end, we subdivide the time interval $[0,T]$ over which we desire to control the system into $K$ subintervals $[t_{k-1},t_k]$ ($1\leq k \leq K$) of equal duration $T/K$,
that is, $t_k = kT/K$ for $0\leq k\leq K$.  In our computational illustrations, we take $K=100$.

If we have evolved the system as far as time $t_k$, we use the following strategy to determine our control over the interval $[t_k,t_{k+1}]$.  As a preliminary step, we solve the deterministic problem on the interval $[t_k,T]$ using the Pontryagin minimum principle. Having determined the control $u(t)$ for this deterministic problem, we apply it to the stochastic problem, but only for the time interval $[t_k,t_{k+1}]$ and thus we evolve the system to time $t_{k+1}$. Performing this process $K$ times, with the time interval used for the deterministic Pontryagin minimum principle calculation progressively shortening, we arrive at the final time $T$.
By repeating this simulation many times, we can evaluate the average performance of this approach on our stochastic problem. We will use the same process to evaluate a dynamic programming approach, the only difference being the way the control is found. } 

\subsection{Framing our problem as a deterministic optimal control problem}
Let us return to our original problem (Equations \ref{equation1}--\ref{equation2}). As outlined above, we need to find a deterministic optimal control problem capturing the uncertainty on both $x(t)$ and $y(t)$. A natural approach to do so is to replace the states $x(t)$ and $y(t)$ by their expected values. Let us consider the changes to the state $y(t)$ first, which is our belief that Equation \ref{equation1} is true. 

\subsubsection{Changes to the state $y(t)$}

For $t \geq 0$ and $\delta t > 0$ we write $\delta x = x(t+\delta t) -x(t)$. From Equations \ref{equation1} and \ref{equation2} we see that the Wiener component to the displacement $\delta x$ is
\begin{align}
\int_t^{t+\delta t} dB_{t '}=\begin{cases} 
\delta x - \int_t^{t+\delta t} u(t')dt' &\text{if Equation \ref{equation1} is true};\\
\delta x + \int_t^{t+\delta t} u(t')dt' & \text{if Equation \ref{equation2} is true} 
\end{cases}
\end{align}
and so the relative likelihoods of the observed state change $\delta x$ are
\begin{align*}
\frac{1}{\sqrt{2\pi\delta t}}\exp\Bigl\{-\frac{1}{2\delta t}\Bigl[ \delta x - \int_t^{t+\delta t} u(t')dt' \Bigr]^2\Bigr\}&~\text{if Equation \ref{equation1} is true};\\
\frac{1}{\sqrt{2\pi\delta t}}\exp\Bigl\{-\frac{1}{2\delta t}\Bigl[ \delta x + \int_t^{t+\delta t} u(t')dt' \Bigr]^2\Bigr\}&~\text{if Equation \ref{equation2} is true}.
\end{align*}
Based on prior beliefs $y(t)$ and $1-y(t)$ for Equations \ref{equation1} and \ref{equation2}, respectively, we can update the information state $y(t+\delta t)$ using Bayes' theorem. The updated information state is calculated as the multiplication of the prior belief $y(t)$ and the relative likelihood of the observed state change $\delta x$ and is then normalized:
\begin{align}
y(t+\delta t) &= 
\frac{\displaystyle \frac{y(t)}{\sqrt{2\pi\delta t}}\exp\Bigl\{-\frac{1}{2\delta t}\Bigl[ \delta x - \int_t^{t+\delta t} u(t')dt' \Bigr]^2\Bigr\}}
{\displaystyle\frac{y(t)}{\sqrt{2\pi\delta t}}\exp\Bigl\{-\frac{1}{2\delta t}\Bigl[ \delta x - \int_t^{t+\delta t} u(t')dt' \Bigr]^2\Bigr\} + \frac{1-y(t)}{\sqrt{2\pi\delta t}}\exp\Bigl\{-\frac{1}{2\delta t}\Bigl[ \delta x + \int_t^{t+\delta t} u(t')dt' \Bigr]^2\Bigr\}} \label{eq:bayes}\\
&=\frac{\displaystyle y(t)\exp\Bigl[ \frac{2\delta x}{\delta t}\int_t^{t+\delta t} u(t')dt' \Bigr]}
{\displaystyle y(t)\exp\Bigl[ \frac{2\delta x}{\delta t}\int_t^{t+\delta t} u(t')dt'\Bigr] + 1 - y(t)}.
\end{align}
It follows from this that
\begin{equation}
y(t+\delta t) - y(t) = \eta\Bigl( y(t),\frac{2\delta x}{\delta t}\int_t^{t+\delta t} u(t')dt' \Bigr), \label{eq:deltaY}
\end{equation}
where for brevity we have written
\begin{equation}
\eta(y,v) = \frac{y(1-y)(e^v - 1)}{ye^v + 1 - y}.
\end{equation}
We note for later use that
\begin{align}   
\eta(y(t),0)&=0,&\frac{\partial\eta}{\partial v}(y(t),0)&=y(t)[1-y(t)],&
\frac{\partial^2\eta}{\partial v^2}(y(t),0)&=y(t)[1-y(t)][1-2y(t)].\label{e:HelpfulDerivatives}
\end{align} 
We assume that $u(t)$ is right-continuous, which implies that 
\begin{equation}
\lim_{\delta t\downarrow 0}\frac{1}{\delta t}\int_t^{t+\delta t} u(t')dt' =u(t).
\end{equation}
In the passage from Equation \ref{eq:deltaY} to a stochastic differential equation, in which $y(t+\delta t) - y(t)$ becomes $dy(t)$, 
we replace the second argument of the function $\eta$ in the right-hand side of Equation \ref{eq:deltaY} by $2u(t) dx(t)$. To be strictly correct this requires us to ask a little more of $u(t)$ than right-continuity. Right differentiability, the slightly weaker requirement of right Lipschitz continuity, or the even weaker requirement of right H\"older continuity with associated exponent $h>1/2$ would suffice \change{(the need for $h>1/2$ is to ensure that the correction to the leading term the integral, when multiplied by $\delta x$, is $o(\delta t)$ and so 
can be ignored in subsequent It\^o calculus)}. We arrive at
\begin{equation}
dy(t)=\eta(y(t),2u(t) dx(t)).\label{e:stochasticDE-before-expansion}
\end{equation}
Further progress requires us to specify which of Equation \ref{equation1} or \ref{equation2} is correct, so we write $dx(t)=\pm u(t)dt+dB_t$ where the upper sign is taken if 
Equation \ref{equation1} is correct and the lower sign is taken if Equation \ref{equation2} is correct. Thus we have
\[
2u(t)dx(t) = \pm 2u(t)^2dt + 2u(t)dB_t,
\] 
which is a diffusion with drift $\mu_t=\pm 2u(t)^2$ and standard deviation $\sigma_t = 2u(t)$. 
Following the normal approach in It\^o calculus we expand $\eta(y(t),2u(t) dx(t))$ to second-order in the second argument using the results \ref{e:HelpfulDerivatives} noted above, replace $(dB_t)^2$ by $dt$ and retain only those terms multiplied by a single factor of $dt$ or $dB_t$, giving
\begin{align}
dy(t)
& = 2y(t)[1-y(t)][1 \pm 1  - 2y(t)]u(t)^2 dt + 2y(t)[1-y(t)]u(t)dB_t.
\end{align}
In order to frame the problem as a deterministic optimal control problem, we now calculate the expected value of this derivative by removing the Wiener process:
\begin{equation}
\mathbf{E} \left[ \frac{dy(t)}{dt} \right] =  \begin{cases}
4y(t)[1 - y(t)]^2u(t)^2&~\text{if Equation \ref{equation1} is true};\\
-4 [1 - y(t)]y(t)^2u(t)^2&~\text{if Equation \ref{equation2} is true}.
\end{cases}
\label{eq:deterdy1and2}
\end{equation}
This result can be interpreted in the terminology of dynamical systems, if we write a deterministic $dy/dt$ in place of the expectation of the random $dy/dt$. The states $y(t)=0$ and $y(t) =1$ (which correspond to certainty about Equation \ref{equation2} being true 
or about Equation \ref{equation1} being true, respectively) are equilibria. However, if $0<y(0)<1$, then if Equation \ref{equation1} is true, we have $dy/dt>0$, and our confidence in the truth of Equation \ref{equation1} increases over time, while if  Equation \ref{equation2} is true, we have $dy/dt<0$, and our confidence in the truth of Equation \ref{equation1} decreases over time, while correspondingly our confidence in the truth of Equation \ref{equation2} increases. Loosely speaking, in an average sense our evolving confidence is attracted to the fixed point of truth.  However, the closer we approach the truth, the more slowly we learn.
For example, if Equation \ref{equation1} is true, then $y(t)\geq y(0)$ for all $t\geq 0$, so we have
\[
\frac{dy(t)}{dt}\geq 4y(0)[1-y(t)]^2u(t)^2
\]
and we can integrate to deduce that
\[
1- y(t) \leq  \frac{1-y(0)}{\displaystyle 1 + 4y(0)[1 - y(0)]\int_{0}^t u(\tau)^2d\tau}.
\]
This shows that to achieve high confidence (that is, $y(t)\approx 1$ if Equation \ref{equation1} is true) from an uncertain initial state, we need very long experience with a weak control  ($|u(t)|\ll 1$), but a shorter time interval for a stronger control.

The analysis of $y(t)$ has led to the deterministic differential equation \ref{eq:deterdy1and2}, which has two different forms depending on which of Equation \ref{equation1} or \ref{equation2} is true.
We follow the natural approach that consists of creating two information states $y_1$ and $y_2$ governed by the corresponding cases that arise in Equation \ref{eq:deterdy1and2}:
\begin{align}   
&\frac{dy_1(t)}{dt} = 4y_1(t)[1-y_1(t)]^2u^2(t),  \\
&\frac{dy_2(t)}{dt} = -4y_2(t)^2[1-y_2(t)]u^2(t).
\end{align}
The states $y_1$ and $y_2$ will be part of the deterministic optimal control problem. The two states have the same initial condition: $y_1(0)=y_2(0)=y(0)$. This concludes our analysis of the state $y(t)$.  

\subsubsection{Changes to the state $x(t)$}

We now show how to deal with the state $x(t)$. First, we argue that we only need to deal with the case where $x(t)$ is positive, because the problem is symmetric. If $x(t) < 0$, the two possible true equations can be written
\begin{align} 
d(-x(t)) &= -(u(t)dt + dB_t)= (-u(t))dt - dB_t, \ \ \ \ \ \ \ \ \ \ (\text{if Equation } \ref{equation1}\text{ is true})    \\
d(-x(t)) &= -(-u(t)dt + dB_t)= -(-u(t))dt - dB_t,  \ \ \ \ \  (\text{if Equation } \ref{equation2}\text{ is true}) 
\end{align}
Because the Wiener process has no drift, $-dB_t$ has the same distribution as $dB_t$. So, the state equation governing $-x(t)$ when applying $-u(t)$ is equal in distribution to the state equation governing $x(t)$ when applying $u(t)$. Also, we can change the sign of $u(t)$ safely because the control space is of the form $[-U,U]$ and the objective depends on $u(t)^2$ and thus does not depend on the sign of $u(t)$. Hence, denoting $u(t)$ an optimal control (there might be many) for the state $x(t)$ and $y(t)$, $-u(t)$ is also an optimal control on states $-x(t)$ and $y(t)$. The optimal control is thus of the form 
\begin{align}
u(t)=\sgn(x(t))f_u(|x(t)|,y(t),t),  \label{eq:formOfU}
\end{align} 
where the function $f_u$ is to be determined. So, we need only consider $z(t):=|x(t)|$ in the deterministic problem.

We want the deterministic differential equation of $z(t)$ to account for three sources of randomness and uncertainty. Firstly, the true equation is `selected' randomly at $t=0$. Secondly, our knowledge of the true equation is imperfect. Thirdly, each equation is stochastic due to the Wiener process. We treat these three aspects in order.

Firstly, the randomness of the true equation is perhaps the simplest to deal with. Since Equations \ref{equation1} and \ref{equation2} are respectively true with probabilities $y({0})$ and $1-y({0})$, the derivative of $z(t)$ will be of the form
\begin{align}
\frac{dz(t)}{dt} = y({0}) \frac{dz_1(t)}{dt} + [1-y({0})] \frac{dz_2(t)}{dt},  \label{eq:dz}
\end{align} 
where $z_1$ and $z_2$ represent $|x(t)|$ when Equation \ref{equation1} or \ref{equation2} is true, respectively. 

Secondly, let us consider our imperfect knowledge of the true equation. Because the true equation is unknown, we define $dz_1(t)/dt$ and $dz_2(t)/dt$ based on our beliefs of the true equation, rather than on the true equation itself. For example, if we assume that Equation \ref{equation1} is true, our beliefs in Equation \ref{equation1} and \ref{equation2} are $y_1(t)$ and $1-y_1(t)$, which yields: 
\begin{align}
 \frac{dz_1(t)}{dt} &:=\left\{ y_1(t) \mathbf{E}\left[\left. \frac{d |x(t)|}{dt}\right\vert \mbox{Equation\ \ref{equation1} true}\right] + [1-y_1(t)] \mathbf{E}\left[\left. \frac{d |x(t)|}{dt}\right\vert \mbox{Equation\ \ref{equation2} true}\right] \right\} \\
&=\sgn(x(t))\left\{ y_1(t) \mathbf{E}\left[\left. \frac{d x(t)}{dt}\right\vert \mbox{Equation\ \ref{equation1} true}\right] + [1-y_1(t)] \mathbf{E}\left[\left. \frac{d x(t)}{dt}\right\vert \mbox{Equation\ \ref{equation2} true}\right] \right\} \\
&= \sgn(x(t)) \bigl[2 y_1(t) - 1\bigr] u(t) \\
&= \bigl[2 y_1(t) - 1\bigr] f_u(|x(t)|,y_1(t),t),
\end{align}
where we have used Equation \ref{eq:formOfU} \change{and have also noted that for $x(t)\neq 0$ we have $d |x(t)| /dt= \sgn(x(t)) dx(t)/dt$.}
Recall that we always want to reduce $|x(t)|$ to minimize costs, so we would like to have $dz_1(t)/dt\leq 0$. We can do this by setting
\begin{align}
\sgn(f_u(|x(t)|,y_1(t),t)) = - \sgn[2 y_1(t) - 1].
\end{align}   
We conjecture that this sign is always optimal, which seems very sensible because there is no incentive to increase $|x(t)|$ instead of decreasing it. The equation becomes
\begin{align}
\frac{dz_1(t)}{dt} = -|2 y_1(t) - 1| |u(t)|.  \label{eq:dz1}
\end{align}   
For $z_2$, we have
\begin{align}
\frac{dz_2(t)}{dt} = -|2 y_2(t) - 1| |u(t)|.   \label{eq:dz2}
\end{align}   
Finally, combining Equations \ref{eq:dz}, \ref{eq:dz1} and \ref{eq:dz2} leads to
\begin{align}
\frac{dz(t)}{dt} &= - |u| \left\{ y({0})\left| 2 y_1(t)-1\right| +[1-y({0})] \left| 2 y_2(t)-1\right|\right\}   . \label{eqU}
\end{align}
Note that $z_1$ and $z_2$ do not appear in the deterministic model---only $z$ does. This model appears to capture our uncertainty about $y_{true}$. At $t=0$, the derivative is $-|u||2 y({0})-1|$, which is small for `uncertain' values of $y({0})$ (around 0.5). This derivative can be seen as `bridled' in order to penalize poor knowledge. When $t$ increases, the beliefs $y_1$ and $y_2$ converge to 0 and 1, so the derivative converges to  $-|u|$. Controls will tend to have a higher impact on $z(t)$ when $t$ increases, which reflects our increasing knowledge of the true equation and thus better control of the system.  

Thirdly, let us address the uncertainty of $x(t)$ due to the Wiener process. The above formula is insufficient: if $x({0})=0$, $z(0)=0$ so the optimal control is zero for a total cost of zero. In reality, $x(t)$ varies according to the Wiener process and the total cost is positive almost surely. To account for this, we denote by $\xi(t) =\mathbf{E} |x(t)|$ the expected modulus of an uncontrolled $x(t)$, with $x(0)=0$. Appropriate integration of the Gaussian probability density function with mean zero and variance $t$ establishes that $\xi(t) = (2t/\pi)^{1/2}$, from which it follows that
\begin{align}   
\frac{d\xi(t)}{dt} = \frac{1}{\pi\xi(t)}.   \label{eqXi}
\end{align}
Although this differential equation is an approximation when $x(0) \neq 0$, we expect it to remain accurate in our case because a successful control will drive the process towards small values of $x$. 

In order to account for all three aforementioned sources of uncertainty, we defined the differential state equation of $z(t)$ in the deterministic problem as the sum of Equations \ref{eqU} and \ref{eqXi}:
\begin{align}   
&\frac{dz(t)}{dt} \approx \frac{1}{\pi z(t)} - u \left\{ y({0})\left| 2 y_1(t)-1\right| +[1-y({0})] \left| 2 y_2(t)-1\right|\right\}. \label{eq:modelDeterm}
\end{align}
Note that combining equations this way is also an approximation, because the first term $1/\pi z(t)$ corresponds to the uncontrolled case. However, we can observe that this differential equation achieves to capture the antagonism between our control under uncertainty (second term), which reduces $z(t)$, and the `penalty' caused by the stochasticity of the Wiener process (first term), which increases $z(t)$. The rest of the deterministic problem is
\begin{align}   
&\frac{dy_1(t)}{dt} = 4y_1(t)[1-y_1(t)]^2u^2(t),  \\
&\frac{dy_2(t)}{dt} = -4y_2(t)^2[1-y_2(t)]u^2(t), \\
&\min_{u} \int_0^T [z^2(t) + u^2(t)] dt, \label{eq:modelDetermEnd} 
\end{align}
with initial conditions 
\begin{align}   \label{eq:InitialConditions}
z(0)=|x({0})|, y_1(0) = y_2(0) = y({0}).
\end{align}  
We restrict $u$ to only non-negative values because the choice of sign has been made already. Interestingly, the trade-off on the control $u$ appears perhaps more clearly in this deterministic problem than in the stochastic one, because $u^2$ is both a cost to minimize and a linear factor in the derivative of $y_1$ and $y_2$. In other words, `extreme' controls are costly but increase our knowledge, potentially leading to better future control. This is a classic trade-off in the dual control and adaptive management literature. 

\subsection{Applying the Pontryagin minimum principle}

We can now apply the Pontryagin minimum principle to our deterministic problem with states $z(t)$, $y_1(t)$ and $y_2(t)$ (more details to be found in Appendix \ref{appendix:pontryagin}). The Hamiltonian is: 
\begin{align}    \label{eq:Hamiltonian}
&H = g + \lambda h, \notag \\
&=u(t)^2+z(t)^2 + \lambda_z(t) \left\{ \frac{1}{\pi z(t)} - u(t) y({0})\left| 2 y_1(t)-1\right| -u(t) [1-y({0})] \left| 2 y_2(t)-1\right|\right\},  \notag \\
&\ \ \ + 4 \lambda_{y_1}(t) y_1(t)[1-y_1(t)]^2u^2(t) -4 \lambda_{y_2}(t) y_2(t)^2[1-y_2(t)]u^2(t). 
\end{align} 
The adjoints $\lambda_z$, $\lambda_{y_1}$ and $\lambda_{y_2}$ satisfy the equations:
\begin{align}  \label{eq:AdjointEquations}
&\frac{d \lambda_z}{dt} = - \frac{\partial H}{\partial z} = \frac{\lambda_z(t)}{\pi z(t)^2} - 2z(t), \\
&\frac{d \lambda_{y_1}}{dt} = - \frac{\partial H}{\partial y_1} = 2 u(t) \lambda_z(t) y({0}) \sgn[2 y_1(t)-1] -4 \lambda_{y_1}(t) [3 y_1(t)^2-4 y_1(t)+1] u(t)^2, \\
&\frac{d \lambda_{y_2}}{dt} = - \frac{\partial H}{\partial y_2} = 2 u(t)  \lambda_z(t) [1-y({0})] \sgn[2 y_2(t) - 1] +4 \lambda_{y_2}(t) y_2(t) [2-3 y_2(t)] u(t)^2 \label{eq:AdjointEquationsEnd},\\
&\lambda_{z}(T)=0,\ \lambda_{y_1}(T)=0,\ \lambda_{y_2}(T)=0. \label{eq:FinalConditions}
\end{align} 
The optimal control satisfies: 
\begin{align}
u(t) &= \argmin_{u \in [0,U]} H(u). \label{eq:minU}
\end{align} 
\change{The optimal control problem has been written out for the initial step, where the process starts at time $0$ and the time interval over which control is sought is $[0,T]$.
As explained in Section \ref{secApproachOverview}, this approach is used on progressively shorter time intervals $[t_k,T]$ as the algorithm proceeds.}

\subsection{Multidimensional formulation}  \label{section:multidim}

We can also consider a higher-dimensional version of this problem, where we have $N$ states and $N$ controls. In this case, we assume that one of the following two sets of stochastic differential equations is true:
\begin{align} 
dx_i(t) &= u_i(t)dt + dB_t \text{   for all } 1 \leq i \leq N, \label{eqmultidim1} \\
dx_i(t) &= -u_i(t)dt + dB_t \text{   for all } 1 \leq i \leq N,   \label{eqmultidim2}
\end{align}
where $B_t$ is a multidimensional Wiener Process. Again we aim to keep the states and controls small in a suitable mean-square sense over the time interval $[0,T]$:
\begin{align}
\lefteqn{\min_{u} \Bigl\{\dfrac{1}{2}\ \mathbf{E} \Big[
\int_0^T \left[ \sum_{i=1}^N u_i^2(t) + \sum_{i=1}^N x_i^2(t) \right] dt \mid \mbox{Equation\ \ref{equation1} true} \Big]}\nonumber\\
&\qquad  + \dfrac{1}{2}\ \mathbf{E} \Big[
\int_0^T \left[ \sum_{i=1}^N u_i^2(t) + \sum_{i=1}^N x_i^2(t) \right] dt \mid \mbox{Equation\ \ref{equation2} true} \Big] \Bigr\}.
\label{equation3multi}
\end{align}
The deterministic problem becomes:
\begin{gather}   
\frac{dz_i(t)}{dt}=\frac{1}{\pi z_i(t)} - u_i \left\{ y({0})\left| 2 y_1(t)-1\right| +[1-y({0})] \left| 2 y_2(t)-1\right|\right\}  \text{   for all } 1 \leq i \leq N,\\
\frac{dy_1(t)}{dt} = 4y_1(t)[1-y_1(t)]^2\sum_{i=1}^N u_i^2(t),   \\
\frac{dy_2(t)}{dt} = -4y_2(t)^2[1-y_2(t)]\sum_{i=1}^N u_i^2(t), \\
\min_{u} \int_0^T \left[ \sum_{i=1}^N u_i^2(t) + \sum_{i=1}^N z_i^2(t) \right] dt.
\end{gather}

\subsection{Obtaining the optimal deterministic control}

We have shown how to cast our dual control problem into a deterministic optimal control problem, which is our main theoretical contribution. Several methods exist to solve deterministic optimal control problems. We use the forward-backwards sweep. 

In short (more details in Appendix \ref{appendix:fbsweep}), this approach takes as arguments the initial states $x({0}), y({0})$ and a time horizon $T'$ (Algorithm\ \ref{algodeterm}). \change{(As the global algorithm proceeds, the initial states become $x(t_k), y(t_k)$ and $T' = T - t_k$.)} Based on an initial guess $u$ for the control function, the forward sweep solves the state equation (Equation \ref{eq:state}, Line \ref{line:odestate} in Algorithm\ \ref{algodeterm}) forwards in time. It takes as arguments the initial states $|x({0})|, y({0}), y({0})$ for the functions $z$, $y_1$ and $y_2$, respectively, and the time horizon $T'$ and the control profile $u$. Then, the backward sweep solves the adjoint equations (Equations \ref{eq:AdjointEquations}--\ref{eq:AdjointEquationsEnd}, Line \ref{line:odeadjoint}) backwards in time, taking as argument the final costates (all equal to 0), the time horizon $T'$, the control profile $u$ and the values of the states $z, y_1, y_2$ throughout the time-frame. Finally, at regular times $t$ in the time-frame, the controls $u(t)$ are updated so as to minimize $H(x(t), u(t), \lambda(t))$. Rather than using a simple update for $u(t)$, we update $u(t)$ with a weighted sum of the new estimate and the previous estimate. The weight parameter $\omega_\tau$, used to improve numerical stability, favours the new estimate early in the computation and the previous estimate late in the calculation. This process is repeated until convergence.

\begin{algorithm}
\caption{SolveDeterministic($x({0}), y({0}), T'$)}   \label{algodeterm}
\begin{algorithmic}[1]
\State \textbf{Initialization}: $\tau = 0, u(t) := 0, u_{c}(t) := 1\mbox{\emph{ for all }}t \in [0,T']$ \label{line:initialcontrol}
\While {$||u-u_{c}||^2 \geq \epsilon$}     \label{line:loop}
\State $z, y_1, y_2 := \textit{ForwardSweep\/}(|x({0})|, y({0}), y({0}), T', u)$   \label{line:odestate}
\State $\lambda_z, \lambda_{y_1}, \lambda_{y_2}  := \textit{BackwardSweep\/}(0, 0, 0, T', u, z, y_1, y_2) $    \label{line:odeadjoint}
\For {$t = 1 :  \textit{nTimeSteps\/}$}
\State $u_{c}(t):= \argmin_{\widetilde{u}}\Bigl[ \widetilde{u}^2 - \lambda_z(t) \bigl( \widetilde{u} y({0})\left| 2 y_1(t)-1\right| -\widetilde{u} [1-y({0})] \left| 2 y_2(t)-1\right|\bigr) + 4 \lambda_{y_1}(t) y_1(t)[1-y_1(t)]^2\widetilde{u}^2 -4 \lambda_{y_2}(t) y_2(t)^2[1-y_2(t)]\widetilde{u}^2 \Bigr]$   \label{line:findU}
\State $u(t):= \omega_\tau \cdot u_{c}(t) + (1- \omega_\tau) \cdot u(t)$   \label{line:updateU}
\EndFor
\State $\tau := \tau + 1$
\EndWhile
\State \Return the control profile $u(t)$
\end{algorithmic}
\end{algorithm}

\subsection{Evaluation and benchmarking}

\subsubsection{Comparison to dynamic programming}

Dynamic programming is a popular solution method for dual control problems, though it is not necessarily the most efficient approach and
for large problems its implementation may not be practical.  When it can be implemented, it provides truly optimal solutions.
We shall use dynamic programming to test how well our continuous-time approach performs.
Additionally, by analyzing computational times and memory requirements, we illustrate how our method avoids the curse of dimensionality, which is the main motivation behind this work.

To use dynamic programming, a few steps are required (see Appendix \ref{appendix:momdpanddynamicprogramming}, in particular Sections \ref{momdpframe} and \ref{appendix:dynamicprogramming}, for full details). The most important steps are:
\begin{itemize}
    \item discretizing the action and state spaces, including the information state;
    \item discretizing the time-frame into $K=100$ regular time intervals. We map the time intervals $[0,t_1], [t_1,t_2], \ldots, [t_{K-1},T]$ to time steps $0,1,\ldots,K-1$ in dynamic programming;
    \item finding the optimal policy by working backwards (Algorithm\ \ref{DynamicProgramming}), using the Principle of Optimality \cite{bellman_dynamic_1957} that all remaining actions must constitute an optimal policy relative to the current state.
\end{itemize}
The main differences with our optimal control approach are:
\begin{itemize}
    \item Dynamic programming is guaranteed to find the global optimum for the given discrete space, but not for the original continuous problem -- this is particularly relevant with coarse discretizations in space and/or time.
    \item Dynamic programming runs offline and returns the entire policy. In contrast, optimal control `only' returns a control profile and runs online -- more on this in the next section.
    \item The main caveat of dynamic programming is the curse of dimensionality, i.e. the number of states handled by dynamic programming increases exponentially with the dimensions of the state space.
\end{itemize}

\begin{algorithm}[h]
\caption{DynamicProgramming($K$)}   \label{DynamicProgramming}
\begin{algorithmic}[1]
\State \textbf{Initialization}: $V(x, y, K) :=  0\mbox{\emph{ for all states }} x \in \widetilde{X},
  y \in \widetilde{Y}$
\For {$t = K-1 : 0$} \ \ \ \ \ \ \ \ \ // find optimal policy backwards
\For {$x \in \widetilde{X},  y \in \widetilde{Y}$}
\For {$a \in \widetilde{A}$}
\State $r(x, y,a) = T/K \cdot (a^2(t) + x^2(t))$ 
\EndFor 
\State $V(x, y, t):=\max\limits_{a \in \widetilde{A}}\  \Bigl[  r(x, y,a) + \sum\limits_{x' \in \widetilde{X},  y' \in \widetilde{Y}}  P(x',  y'|x,  y) V(x',  y', t + 1)\Bigr] $
\State $\pi(x, y, t):=\argmax\limits_{a \in \widetilde{A}}\  \Bigl[  r(x, y,a) + \sum\limits_{x' \in \widetilde{X},y' \in \widetilde{Y}}  P(x',  y'|x,  y) V(x',  y', t + 1)\Bigr] $
\EndFor 
\EndFor 
\State \Return the policy $\pi$
\end{algorithmic}
\end{algorithm}

\subsubsection{General algorithm}

Both the dynamic programming and optimal control approaches are evaluated through simulations with a shrinking horizon (Algorithm \ref{simulations}). 
\change{We recall that we subdivide the time interval $[0,T]$ over which we desire to control the system into $K$ subintervals $[t_{k-1},t_k]$ ($1\leq k \leq K$) of equal duration $T/K$,
that is, $t_k = kT/K$ for $0\leq k\leq K$.  In our computational illustrations, we take $K=100$.}
At the beginning of each simulation, the true state equation, i.e. the variable $y_{true}$, is drawn randomly (Line\ \ref{line:draw}). The control is obtained from dynamic programming or from the optimal control approach (Line\ \ref{line:SolveDeterministic}), depending on which approach is being evaluated. There is a key difference between these two approaches:
\begin{itemize}
\item Since dynamic programming runs offline, the control is simply equal to the policy $\pi$ corresponding to the nearest discretized $\widetilde{x}$ and $\widetilde{y}$.
\item In contrast, our optimal control approach (Algorithm\ \ref{algodeterm}) runs online: it is called at each time step of each simulation with the current values of $|x(t)|$ and $y(t)$, and outputs the control to implement.
\end{itemize}
The next state is drawn based on this control and on the true state equation (Equation \ref{equation1} or \ref{equation2}, Line\ \ref{line:drawState}); the variable $y(t)$ is updated through Bayes' theorem (Equation \ref{eq:bayes}, Line\ \ref{line:updateQ}). The time step is updated until the end of the time horizon (Line\ \ref{line:for}). The output is the average cost per simulation.
\begin{algorithm}
\caption{SolveStochastic($x(0), y(0), T, K$)}   \label{simulations}
\begin{algorithmic}[1]
\For {$i = 1 : nSimulations$}
\State \textbf{Initialization}: $Cost := 0$
\State $y_{true} := DrawQ(y(0), 1-y(0))$   \label{line:draw}
\For {$k = 0 : K-1$}  \label{line:for}
\State $u(t_k) := \pi(x(t_k), y(t_k), k)$ or $u(t_k) := \emph{SolveDeterministic}(x(t_k), y(t_k), T-t_k)$ \ \ \ // from dyn. prog. or opt. control \label{line:SolveDeterministic}
\State $x(t_{k+1}) := \emph{DrawNextState}(u(t_k), x(t_k), y_{true})$  \label{line:drawState} 
\State $y(t_{k+1}) := \emph{UpdateKnowledge}(y(t_k), x(t_k), x(t_{k+1}))$  \label{line:updateQ}
\State $Cost := Cost + T/K \cdot (u^2(t_k) + x^2(t_k))$
\EndFor 
\EndFor 
\State \Return the average cost $Cost/nSimulations$
\end{algorithmic}
\end{algorithm}

\subsubsection{Computational experiments} \label{computational_experiments}

We run simulations for both the dynamic programming and the optimal control approaches, with various initial values for $x(0)$ and $y(0)$ (Table \ref{tab:results}). We show the average cost, 95\% confidence interval and computational time over 500 simulations for each instances.

Note that comparing dynamic programming to the optimal control approach is not competely straightforward, because dynamic programming is offline while the optimal control approach is online. Dynamic programming calculates the optimal policy for all discretized states before the simulations start (offline). So, the simulations are very fast because decisions can be found in a lookup-table. In contrast, our optimal control approach calculates the best control during the simulations by solving an optimal control problem (online). There is no preprocessing time but the simulations are slower than for dynamic programming. The times shown in Table \ref{tab:results} correspond to the sum of the preprocessing and simulation times.

\subsubsection{Results}

\change{For the reason exposed in Section \ref{computational_experiments}, the computation times for both approaches are difficult to compare, so we shall only compare the costs produced by their selected controls. For 1- or 2-dimensional problems, both approaches yield similar costs, i.e. within one another's 90\% confidence intervals.} With 3 or more dimensions, dynamic programming is either intractable or yields significantly higher costs than optimal control due to a too coarse discretization (43.0 vs 35.3 with 3 dimensions). 

Further, the optimal control approach appears to do better in higher dimensions: the computational time per dimension decreases when the dimension increases (e.g. for dimensions 1 and 10, 39,600s/1 > 216,000s/10), as well as performance costs per dimension (16.7/1 > 124.4/10). This is because the state $y_{true}$ is learned faster with more dimensions, reducing the average cost per dimension. A positive consequence is the decrease in computational time because a better knowledge means the problem is simpler. This is further supported by the problems starting with $y({0}) = 1$, which are solved faster by the optimal control approach than when $y({0}) = 0.5$.

\begin{table}
\centering \small
\begin{tabular}{| l | c | c | c | c | c | } 
\hline  Approach & \multicolumn{3}{c|}{Dynamic programming (offline)} & \multicolumn{2}{c|}{Optimal control (online)} \\ \hline
Problem & \parbox[t]{2cm}{\centering Cost $\pm$ 95\% \\  confidence} & Time & \parbox[t]{2.2cm}{\centering Discretization \\  $X$/$Y$/$A$} & \parbox[t]{2cm}{\centering Cost $\pm$ 95\% \\  confidence} & \parbox[t]{2cm}{\centering Time (500 \\ simulations)} \\ \hline
  $x({0})=0$ & 
16.7 $\pm$ 1.5 & 43s & 70/45/35 & 16.7 $\pm$ 2.4 & 39,600s \\ \hline
  $x({0})=5$ & 
93.1 $\pm$ 6.0 & 17s & 50/20/15  & 96.6 $\pm$ 7.9 & 43,200s   \\ \hline
  $x({0})=15$ & 
1294.3 $\pm$ 37.0 & 15s & 50/20/15 & 1295.1 $\pm$ 38.4 & 43,200s  \\ \hline
  $x({0})=0$, $U=5$ & 
13.1 $\pm$ 0.6 & 32s & 70/45/35 & 12.8 $\pm$ 0.5 & 45,500s  \\ \hline
  $x({0})=0$, $y({0}) = 1$ & 
11.2 $\pm$ 0.5 & 44s & 70/45/35 & 11.5 $\pm$ 0.7 & 30,000s  \\ \hline
  $x({0})=5$, $y({0}) = 1$ & 
63.6 $\pm$ 3.2 & 17s & 50/20/15 & 60.9 $\pm$ 3.1 & 32,400s  \\ \hline
  $x({0})=(0,0)$ & 
\parbox[t]{2cm}{\centering 99.4  $\pm$ 4.7  (*) \\ 28.2 $\pm$ 1.4 (*)} & \parbox[t]{0.7cm}{\centering  60s \\ 240s} &\parbox[t]{1.5cm}{\centering 10/20/5 \\ 30/20/5 } & 26.6 $\pm$  1.4 & 62,100s  \\ \hline
$x({0})=(0, 0, 0)$ & 
\parbox[t]{3cm} {\centering 43.0 $\pm$ 1.5 (*)\\ Out of memory (*)} &  \parbox{0.7cm}{\centering 420s \\ } & \parbox[t]{1.5cm} {\centering   20/15/5    \\      30/20/5}   & 35.3 $\pm$ 1.3 &  76,200s  \\ \hline
$x({0})=(0, 0, 0, 0)$ & 
Out of memory &  & 20/15/5   & 46.4 $\pm$ 1.3 & 92,700s  \\ \hline
\parbox[t]{3cm} {$x({0})=(0, 0, 0, 0, 0,$  \\ $\mathstrut\qquad 0, 0, 0, 0, 0)$ }& 
Out of memory & & 20/15/5 & 124.4 $\pm$ 2.6 & 216,000s  \\ \hline
\end{tabular} 
\caption{ Average cost, 95\% confidence interval and computational times of both dynamic programming and optimal control, with $y({0}) = 0.5$ and $U=1$ unless otherwise stated. For dynamic programming, the number of discretized bins of $X$, $Y$ and $A$ are shown (number per dimension when $X$ is multidimensional).
Asterisks (*) denote cases where we ran dynamic programming with different sets of discretization parameters $X/Y/A$ to trade off tractability against performance. The memory is set to 10GB.} \label{tab:results}
\end{table} 

\begin{figure}[t]
\centering
\includegraphics[width=\textwidth]{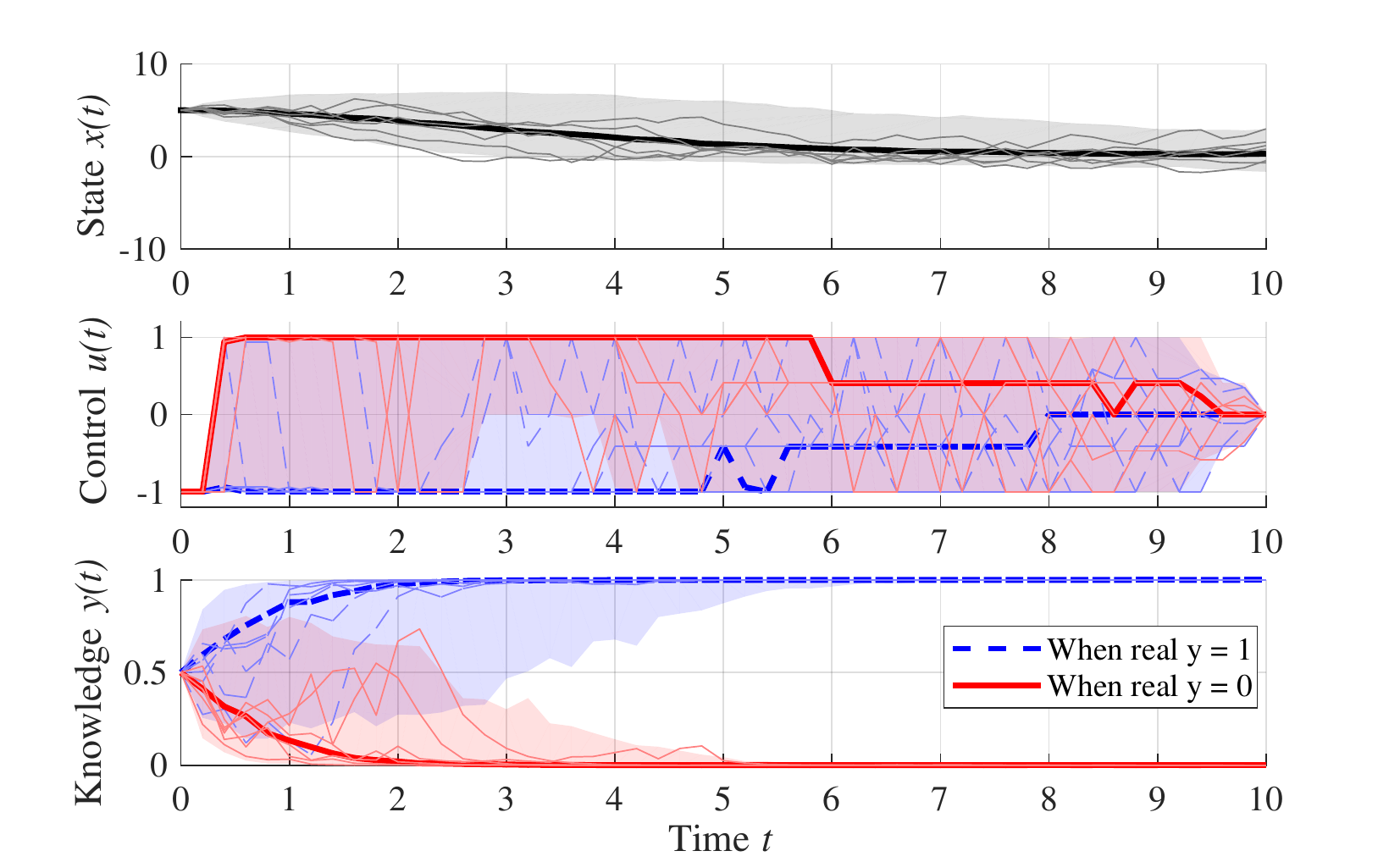} 
\caption{Simulation case study showing the state $x(t)$, the control $u(t)$ and the knowledge $y(t)$ for $0\leq t\leq T=10$, with the constraint $|u(t)|\leq U=1$, initial state $x({0})=5$ and initial knowledge $y({0})=0.5$. We ran 500 simulations. Concerning the control and knowledge, all curves and areas correspond to those of the 500 simulations for which the draw of the variable $y_{true}$ (Line\ \ref{line:draw} in Algorithm \ref{simulations}) is 0, i.e.\ Equation \ref{equation1} is true. The dashed blue line corresponds to simulations for which $y_{true}=1$. For the state $x(t)$, we do not discriminate between $y_{true}=0$ and $y_{true}=1$ because the variations of $x(t)$ do not depend on the value of $y_{true}$ with the initial condition $y({0})=0.5$. In all three graphs, thick curves represent medians, shaded areas show simulations between the $5^{th}$ and $95^{th}$ percentiles, and thin curves show randomly selected individual simulations. The control manages to reduce $x$ until it reaches zero at the end of the time-frame. For $t \leq 2$, the optimal strategy is aggressive (`bang-bang') and selects extreme controls to improve the knowledge. Until $t \approx 6$, the improved knowledge yields extreme controls, causing a steep decrease in the state $x$. Finally, for $t \geq 6$, the state becomes small and controls gradually decreases due to a lack of incentive.}
\label{fig:simuDisplay}
\end{figure}

\begin{figure}[t]
\includegraphics[width=\textwidth]{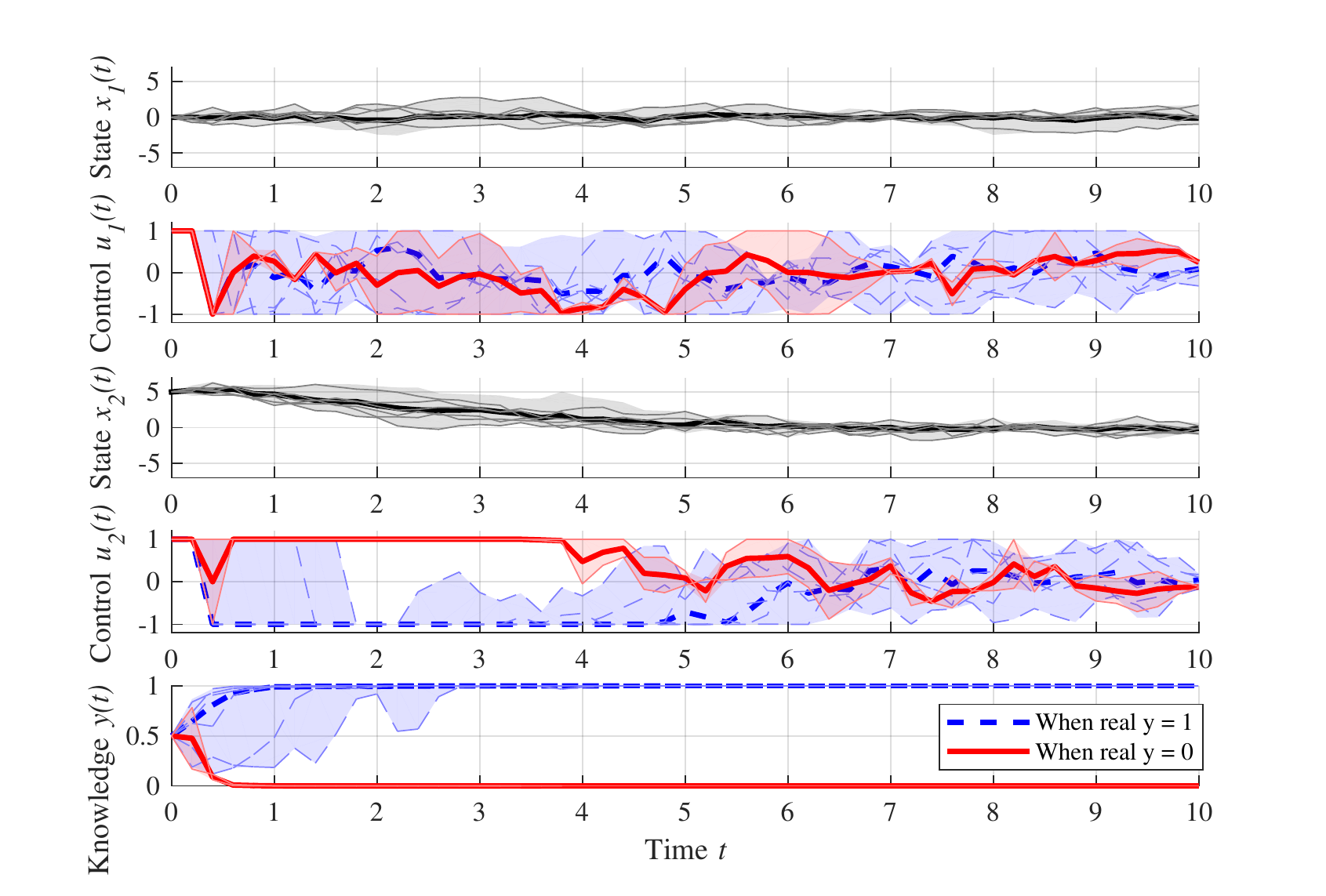} 
\caption{Simulation case study showing the states $x_1(t)$ and $x_2(t)$, the controls $u_1(t)$ and $u_2(t)$ and the knowledge $y(t)$ for $0\leq t\leq T=10$, with the constraint $|u(t)|\leq U=1$ and initial knowledge $y({0})=0.5$. The initial state is $x({0})=(0,5)$, which means $x_1(0)=0$ (top graph) and $x_2(0)=5$ (central graph). Both controls $u_1$ and $u_2$ are extreme at $t=0$, which denotes an aggressive learning strategy. This allows for a quick convergence of the information state $y_{true}$, within roughly one unit of time. Note that as long as $y(t)=0.5$, the expected impact of any control $u_1(t)$ and $u_2(t)$ on $x_1$ and $x_2$ is zero. Reversely, as long as $u_1(t)=0$ and $u_2(t)=0$, $y(t)$ remains $0.5$. So, without an intention to actively learn the information state (i.e. with a passive-adaptive approach, or non-dual in Feldbaum's terminology), we would have no control and no learning throughout the entire time horizon.} \label{fig:simuDisplay2}
\end{figure}

Figures\ \ref{fig:simuDisplay} and \ref{fig:simuDisplay2} show simulations starting with $x({0})=5$ and $x({0}) = (0, 5)$ respectively. The controls are found through dynamic programming but the optimal control approach selects similar levels of control. The solver selects extreme controls in the first time steps despite the uncertainty lying on their consequences, in order to learn quickly. In this simulation the belief converges to the true $y_{true}$ after roughly 2 units of time. With two states, only 1 unit of time is needed to learn the value of $y_{true}$ with quasi-certainty. 

\clearpage

\clearpage

\section{Discussion} 

In this manuscript, we address a continuous-time dual control problem. Given the limitations of dynamic programming, we propose an approach based on optimal control, where the unknown parameter is shown to follow a differential equation. All states are replaced by their expected values, which leads to a deterministic model that is solved with an optimal control algorithm. We evaluate this control profile through simulations in the real stochastic problem. This algorithm rivals dynamic programming on small problems and remains tractable on larger problems, as opposed to dynamic programming. It achieves the right balance between aggressive and smoothly varying controls.

Both the dynamic programming and optimal control approaches have advantages and drawbacks. Dynamic programming is naturally adapted to stochastic problems and handles nonlinearities well. It is also guaranteed to find the optimal solution of discrete-time, discrete-state problems (if tractable), since it explores the entire state space. It is an offline algorithm: it comes with a potentially long preprocessing time, but generates a complete policy. Hence, it suffers from large state spaces, for example when the state space is continuous (therefore infinite) and/or unbounded. It quickly becomes intractable with a multidimensional state (curse of dimensionality \cite{bellman_dynamic_1957}). In these cases, more modeling effort is required to trade off the quality of the solution against the tractability and computational time of the solver. 

In contrast, the optimal control-based approach we propose runs online. This implies that it requires no preprocessing stage. Instead, the algorithm is called at every time step of a simulation at a small computational cost. Additionally, running online means the algorithm has no issues with large and complex state spaces, including multidimensional, unbounded and continuous state spaces, which is a considerable advantage compared to dynamic programming. However, several limitations must be acknowledged. Continuous-time optimal control approaches are not naturally adapted to stochastic problems. Framing the problem as a meaningful deterministic problem required significantly more modeling effort than is needed to run dynamic programming\change{, but unlike model predictive control, we do have the ability to incorporate prior information about system dynamics. 
A possible drawback of our approach is} the lack of numerical performance guarantee, as the solution might only be a local optimum.

The ability to model continuous time and state spaces seems well suited to many real-world problems \change{that
have nearly continuous-time data flow and require very frequent corrective decisions, 
including flight trajectory \cite{kang_pseudospectral_2007}, stock markets \cite{dias_clustering_2015}
(where there may be hidden parameters and models are typically based on Wiener processes),  or medical sciences \cite{hauskrecht_planning_1997} 
(where there are costs to the patient, both from delays in applying therapies while acquiring data to determine the appropriate therapeutic response, and from the commencement of the wrong therapy, that require real-time balancing)}. 
For such problems, continuous-time modeling tools might be more appropriate than discrete-time ones. \change{Not needing} to discretize states \change{is a benefit of our approach.} With a few exceptions, such as presence/absence models, real-world problems usually have continuous state spaces or at least, very large discrete state spaces that need to be partitioned if one is to use dynamic programming; optimal control methods have no such limitations. Further, we can model abundances across many spatial locations without causing dimensionality issues. This can be done by considering either a meta-population \cite{salinas_control_2005}, or by using a spatially explicit optimal control formulation \cite{baker_target_2016,kelly_optimal_2016}. 

However, our approach is too specific to our simple problem to be applied to real-world problems, and would, in that regard, greatly benefit from the following improvements. Firstly, framing the problem as a meaningful deterministic problem is not straightforward. There might be a way to do this more systematically, or to find a stochastic version of the Pontryagin minimum principle that can be directly applied to our problem. Secondly, we exploit the symmetry of our problem with respect to zero. A first step would be to investigate the changes required to the equations and the approach to solve a non-symmetric problem. Thirdly, the uncertainty allowed in our problem is limited to just two possible options. It would be beneficial to handle larger uncertainty sets (one of several equations would be true), or perhaps more realistically, a continuous uncertainty about a parameter, \change{such as Poisson processes, where the frequency of a recurrent event might be unknown, or} $dx(t) = \alpha u(t)dt + dB_t$, with $\alpha \in [-1, 1]$ to be determined. \change{Those types of dual control problems are quite common in other fields, which means our approach could be more easily benchmarked against other methods. For example, in chemical engineering contexts, initial parameter uncertainty within assumed bounds is progressively contracted as time proceeds by nonlinear model predictive control using scenario-tree-based methods or neural nets \cite{thangavel_towards_2015,THANGAVEL201839,LUCIA20141904,DAOSUD20191261,PaulenFIkar2019}. How our approach compares to those methods, in terms of both costs and computation time, would be an interesting question for future research. }

\section*{Acknowledgments}
Christopher Baker is the recipient of a John Stocker Postdoctoral Fellowship from the Science and Industry Endowment Fund. We thank Dan Pagendam for valuable feedback on this manuscript. Computational resources and services used in this work were provided by the HPC and Research Support Group, Queensland University of Technology, Brisbane, Australia.

\bibliography{Peron_et_al_preprint}

\appendix

\section{Solving the problem with Optimal Control}
\label{appendix:optimalcontrol}

\subsection{The Pontryagin minimum principle} \label{appendix:pontryagin}

The Pontryagin minimum principle was developed in the 1960s to solve a class of problems that includes the problem specified by Equations \ref{eq:OCobj1}--\ref{eq:OCobj2}. To this end, one introduces the Hamiltonian function
\begin{align} 
H(x(t), u(t), \lambda(t)) = g(x(t), u(t)) + \lambda(t) h(x(t), u(t)), \label{eq:Hamiltonian_theoretical} 
\end{align} 
where the \textit{adjoint} function $\lambda$, which  plays a similar role to the Lagrange multiplier in simpler optimization problems, satisfies the differential equation
\begin{align} 
\frac{d \lambda}{dt}(t) = - \frac{\partial H}{\partial x} (x(t), u(t), \lambda(t)) \label{eq:adjoint} 
\end{align} 
and the boundary condition
\begin{align} 
\lambda(T) = \frac{d f}{dx}(x(T)),\label{eq:transversality} 
\end{align}  
known as the \textit{transversality condition.}

If there are $n\geq 1$ states, say, we introduce a state vector $\mathbf{x}(t)\in\mathbb{R}^n$ and replace the state equation (Equation \ref{eq:state}) by 
\[
\frac{d \mathbf{x}}{dt}(t) = \mathbf{h}(\mathbf{x}(t),u(t)),
\] 
where $\mathbf{h}:\mathbb{R}^n\times [-U,U]\to\mathbb{R}^n$. In our objective function \ref{eq:OCobj1}, we now have $g:\mathbb{R}^n\times [-U,U]\to\mathbb{R}$ and $f:\mathbb{R}^n\to\mathbb{R}$. Our Lagrange multiplier function becomes $\boldsymbol{\lambda}:[0,T]\to\mathbb{R}^n$ and Equations \ref{eq:Hamiltonian_theoretical}--\ref{eq:transversality} are replaced, respectively, by\[
H(\mathbf{x}(t),u(t),\boldsymbol{\lambda}(t)) = g(\mathbf{x}(t),u(t)) + \boldsymbol{\lambda}(t)\cdot\mathbf{h}(\mathbf{x}(t),u(t)),
\]
\[
\frac{d \boldsymbol{\lambda}}{dt}(t) = -\nabla H(\mathbf{x}(t),u(t),\boldsymbol{\lambda}(t))
\]
 and $\boldsymbol{\lambda}(T) = \nabla f(\mathbf{x}(T))$. Here the gradient operator $\nabla$ acts on the state variable.

The optimal solution is found by minimizing the Hamiltonian with respect to the control and subject to $u(t)$ being feasible (i.e.\ $u(t) \in [-U,U]$ in our case). The problem becomes ill-posed if there is no local minimum of the Hamiltonian in the feasible control interval. Solving Equations \ref{eq:state} and \ref{eq:adjoint} and concurrently minimizing the Hamiltonian yields the optimal control and the corresponding state. If there are $m$ control variables, corresponding to a control vector $\mathbf{u}(t)\in U\subset\mathbb{R}^m$, rather than a single control variable $u(t)$, the procedure is exactly the same.

Since it is often difficult to solve all these equations analytically, numerical methods are standard. One common way is known as the forward-backwards sweep, and is similar to fixed point iteration \cite{hackbusch_numerical_1978}. We start with an initial guess for the control function, and then solve the state equation (Equation \ref{eq:state}) forwards in time. We can then solve the adjoint equation (Equation \ref{eq:adjoint}) backwards in time, using our guess for the control function and the corresponding state. Finally, at regular times $t$ in the time-frame, the controls $u(t)$ are updated so as to minimize $H(x(t), u(t), \lambda(t))$. In our case, minimizing $H$ given $x(t)$ and $\lambda(t)$ is straightforward because $H$ is a one-dimensional quadratic polynomial with respect to $u(t)$---it is either linear, concave or convex. This process is repeated until convergence.

\subsection{Obtaining the optimal deterministic control using the forward-backwards sweep}  \label{appendix:fbsweep}

We now outline how the deterministic control is obtained (Algorithm\ \ref{algodeterm}). Given the current control profile, the trajectories of the states $z$, $y_1$ and $y_2$ are obtained by solving the ordinary differential equation (Equations \ref{eq:modelDeterm}--\ref{eq:modelDetermEnd}, initial conditions in Equation \ref{eq:InitialConditions}) in Line\ \ref{line:odestate}. Then, the adjoints $\lambda_{z}, \lambda_{y_1}$ and $\lambda_{y_2}$ are also found by solving ordinary differential equation (Equations \ref{eq:AdjointEquations}--\ref{eq:AdjointEquationsEnd}, final conditions in Equation \ref{eq:FinalConditions}) in Line\ \ref{line:odeadjoint}. Note that all states, adjoints and controls returned as ODE solutions are discretized in $K=100$ equidistant time steps, but correspond to an exact solution of the ordinary differential equations (by further under-the-hood discretization by the ODE solvers until satisfactory convergence). Based on our experiments, an increase in the number of time steps causes an increase in computational time less than linear.

For given values of the states and co-states, the `candidate' controls $u_{c}$ at each time step is the one which minimizes the Hamiltonian at this time step (Line\ \ref{line:findU} and Eqs.\ \ref{eq:Hamiltonian} and \ref{eq:minU}). However, taken as a whole, the sequence of controls $u_{c}$ rarely equals the \textit{sequential} optimal control because the adjustment can essentially overshoot the optimal point. A more stable update consists of a linear combination between $u_{c}$ with a weight $\omega_\tau$ and the previous control $u$ with a weight $1-\omega_\tau$ (Line\ \ref{line:updateU}), where $\tau$ is the iteration number. We chose for $\omega_\tau$ an exponential decay to allow for significantly changes at the beginning of the algorithm but avoids overshoot after a few iterations. We set $\omega_\tau = e^{-0.15\tau}$ as it achieved a stable and fast convergence over the various instances on which we evaluated this algorithm. The process repeats until the previous and new control are close enough (Line\ \ref{line:loop}). We set the associated threshold $\epsilon$ to $10^{-5}$ as it achieved a good trade-off between accuracy and computational time. The output is the control profile for all time steps. The process can be sped up by setting the initial control profile (Line\ \ref{line:initialcontrol}) as the output control profile from the previous time step in the simulation. 

\section{Solving the problem with Mixed observability Markov decision processes and dynamic programming}
\label{appendix:momdpanddynamicprogramming}
\subsection{Mixed observability Markov decision processes} 

\begin{figure}
\centering
\includegraphics[width=2.5in]{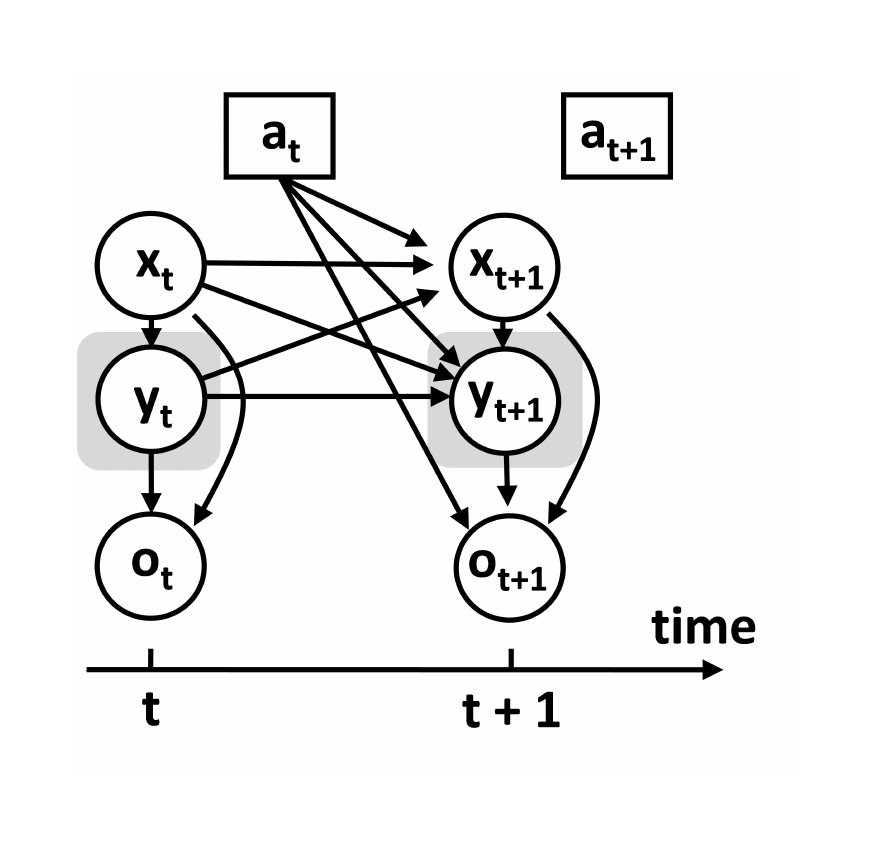}
\caption{Illustration of the interdependencies between states, observations and actions in a MOMDP. The gray area surrounding the variable $y$ indicates that it is partially observed.}
\label{fig:momdplatex}
\end{figure}

A partially observable Markov decision process (POMDP) is a mathematical framework to optimize sequential decisions on a probabilistic system under imperfect observation of the states \cite{sigaud_markov_2010}. Mixed observability Markov decision processes (MOMDPs) are a special case of POMDPs, where the state can be decomposed into a fully observable component, $x$, and a partially observable component, $y$ \cite{ong_planning_2010}. MOMDPs can model various decision problems where an agent knows its position but evolves in a partially observable environment, or when the transition functions or rewards are uncertain. Formally, a MOMDP \cite{ong_planning_2010} is a tuple $\langle X, Y, A, O, P_x, P_y, Z, R\rangle$ with the following attributes.
\begin{itemize}
\item The \emph{state space} is of the form $X \times Y$, with both $X$ and $Y$ of finite cardinality. The current state $(x, y)$ fully specifies the system at every time step. The component $x \in X$ is assumed fully observable and $y \in Y$ is partially observable.
\item The \emph{action space} $A$ is finite.
\item Transition probabilities between states are expressed succinctly using the notational conventions that $(x, y)$ denotes the state immediately before action $a$ is implemented, $(x', y')$ denotes the state immediately after action $a$ is implemented. Also, where the state $y$ is undetermined after the action $a$, it is replaced by an arbitrary variable $v$. We define
\begin{align*}
P_x(x,y,a,x') &= \Pr\Bigl\{\begin{array}{c}(x,y) \mapsto (x',v)~\text{when}\\
\text{$a$ is implemented}\end{array}\Bigr\},\\
P_y(x,y,a,x',y') &= \Pr\Bigl\{\begin{array}{c}(x,y) \mapsto (x',y')~\text{when}\\
\text{$a$ is implemented}\end{array}\,\Bigr\vert\,\begin{array}{c}(x,y) \mapsto (x',v)~\text{when}\\
\text{$a$ is implemented}\end{array}\Bigr\}.
\end{align*}
The process satisfies the Markov property in that these probabilities do not depend on past states or actions.
\item The \emph{reward matrix} is the immediate reward $r(x,y, a)$ that the policy-maker receives for implementing $a$ in state $(x,y)$.
\item The \emph{observation space} $O$ is finite.
\item The \emph{observation probability} $Z$ is defined as
\begin{align}
Z(a,x',y',o') =\Pr\Bigl\{\text{observe}~o'\in O\,\Bigl\vert \begin{tabular}{@{~}c@{~}}
state is $(x',y')$\\
after action $a$
\end{tabular} \Bigr\}.
\end{align}
\end{itemize}
The sequential decision making process unfolds as follows (Figure \ref{fig:momdplatex}). Starting at time $t=0$ in a given initial state $(x({0}), y({0}))$, the decision maker chooses an action $a_0$ and receives the reward $r(x({0}), y({0}), a_0) $. The states $x_1$ and $y_1$ corresponding to $t = 1$ are drawn according to the probabilities $P_x(x({0}), y({0}), a_0, \cdot )$ and $P_y(x({0}), y({0}), a_0, x_1, \cdot )$. The observation $o_1$ is drawn according to the probability $Z(a_0, x_1, y_1, \cdot )$. The decision maker then observes $x_1$ and $o_1$, selects a new action $a_1$ and the process repeats.

The goal of a decision maker is to select actions sequentially to achieve the best expected sum of rewards over time, with respect to a specific optimization criterion. Here, we use a finite time horizon with $K$ time steps, and we seek choices of the actions $a_0,a_1,\ldots, a_{K-1}$ to maximize
\begin{align}
&\mathbf{E}\Bigl[\sum_{t=0}^{K-1}  r(x_t,y_t,a_t)\,\Bigl\vert\, x({0}),y({0})\Bigr].
\end{align}
In the general case, the choice of action $a_t$ may depend on the entire history of actions and observations up to time step $t$ \cite{sigaud_markov_2010,ong_planning_2010}. We will present the concept of a \textit{policy} in the next section more specifically for our context. 

\goodbreak

\subsection{Casting our problem into a discretized MOMDP}
\label{momdpframe}

In order to frame our problem (Equations \ref{equation1}--\ref{equation3}) as a MOMDP, we take the following steps: 
\begin{itemize}
\item We discretize the state $x \in \mathbb{R}$ into a set $\widetilde{X}$, made of regular intervals of size $\delta_x$. Concerning the bounds of $\widetilde{X}$, we can only guarantee that $x$ will stay within a given interval with a certain probability, because the state equations contain a Wiener process $dB_t$. The maximum standard deviation of the Wiener process is achieved at the stopping time $T$ and equals $T^{1/2}$. We set the bounds to three standard deviations for the Wiener process around the initial state $x({0})$, capturing 99.7\% of the values (if the control remains zero). Thus, $\tilde{X}$ is to be confined to the interval $[x({0})-3T^{1/2},x({0})+3T^{1/2}]$. Values over the bounds will be projected back to the bounds. 
\item We discretize the time-frame into $K=100$ regular (equidistant) time intervals, and we denote by $\delta_K=T/K$ the length of each interval. We map these time intervals $[0,\delta_K], [\delta_K,2\delta_K], \ldots, [T-\delta_K, T]$ to time steps $0,1,\ldots,K-1$ in the MOMDP. The MOMDP transition probabilities of each action $a_t$ (time step $t$) are calculated by assuming that, in the continuous-time control problem, the control on the entire interval $[t\delta_K,(t+1)\delta_K]$ is constant and equal to $a_t$. Future states are interpolated on the discretized set $\widetilde{X}$.
\item We discretize the control space into a set $\widetilde{A}$, made of regular intervals of size $\delta_a$. We use the term \textit{action} for these discretized controls, as in the MOMDP literature. 
\item MOMDP solvers are well suited to handling the information state $y$ by themselves: they only require prior belief $y({0})$ on the true value as input (a uniform distribution is standard if with no prior information), without having to manually discretize the belief state. However, the best MOMDP (or, more generally, POMDP) solvers \cite{silver_monte-carlo_2010,kurniawati_sarsop_2008,ong_planning_2010} are tailored for infinite time horizons and our problem has a finite time horizon; we choose instead to apply dynamic programming on a discretized MOMDP. This comes at the cost of discretizing the different states $y \in [0,1]$. The resulting set $\widetilde{Y}$ is of the form $[0,\delta_y, \ldots, 1 - \delta_y, 1]$.
\end{itemize}

\subsection{Dynamic programming} 
\label{appendix:dynamicprogramming}

In dynamic programming, we are looking for an optimal policy $\pi: \widetilde{X} \times \widetilde{Y} \times \{0,1,\ldots,K-1\} \rightarrow \widetilde{A}$. It is a mapping from the discretized time and state spaces to the set of actions (or control) and maximizes the objective criterion. 
To do so, we evaluate the optimal value function $V:\tilde{X}\times\tilde{Y}\times \{0,1,\ldots,K-1\}\to\mathbb{R}$, defined such that $V(x,y,t)$ is the optimal expected sum of rewards received when the system evolves from $(x,y)$ at time $t$ to whatever the final state is at time step $K-1$. Our goal is to calculate $V(x,y,t)$ and the associated optimal policy $\pi(x,y,t)$ for all states $(x,y) \in \widetilde{X} \times \widetilde{Y}$ and time $t\in \{0,1,\ldots,K-1\}$.
We can do this by working backwards (Algorithm\ \ref{DynamicProgramming}), using the Principle of Optimality \cite{bellman_dynamic_1957} that all remaining decisions must constitute an optimal policy relative to the current state. Thus

\begin{align}
V(x, y, K-1) &= \max_{a \in A} r(x,y,a), \\
\pi(x, y, K-1) &= \argmax_{a \in A} r(x,y,a),
\end{align}
and for $t<K-1$, where $P(x',y'|x,y,a)$ is the single step transition probability from state $(x,y)$ to state $(x',y')$ after implementing action $a$, Bellman's equation is satisfied:

\begin{align}
V(x,y,t) &= \max_{a\in A}\Bigl[ r(x,y,a) + \sum_{x',y'\in \tilde{X}\times\tilde{Y}} P(x',y'|x,y,a)V(x',y',t+1)\Bigr], \label{eq:value}\\
\pi(x,y,t)&=\arg \max_{a\in A} \Bigl[ r(x,y,a) + \sum_{x',y'\in \tilde{X}\times\tilde{Y}} P(x',y'|x,y,a)V(x',y',t+1)\Bigr].
\end{align}
The initial action, as prescribed by dynamic programming is $\pi(x({0}),y({0}),0)$, where $y({0})$ is the prior belief at $t=0$.

Note that this policy is optimal in the discretized problem only. It is not guaranteed to be optimal in the real, continuous problem. Its `real' value will be assessed by simulations and will likely differ from the value predicted by Equation \ref{eq:value}. 

Although this dynamic programming approach can work efficiently for small state spaces, it suffers from large or multidimensional state spaces (the curse of dimensionality \cite{bellman_dynamic_1957}). This motivates us to find an alternative approach using tools from optimal control theory.

\end{document}